\documentclass[11pt]{article}

\usepackage{latexsym,graphpap,amssymb,graphics}

\makeatletter
 
\@addtoreset{equation}{section}

\makeatother

\newtheorem{theorem}[equation]{Theorem}
\newtheorem{corollary}[equation]{Corollary}
\newtheorem{definition}[equation]{Definition}

\newtheorem{lemma}[equation]{Lemma}
\newtheorem{proposition}[equation]{Proposition}
\newtheorem{remark}[equation]{Remark}
\newenvironment{rem}{\begin{remark}\rm }{\end{remark}}

\newcommand{\gee}{\mathcal{G}}
\newcommand{\la}{\Lambda}
\newcommand{\eps}{\epsilon}

\newcommand{\lat}{\ensuremath{\mathbb{L}}}
\newcommand{\zed}{\ensuremath{\mathbb{Z}^d}}
\def\reff#1{(\ref{#1})}

\begin{document}

\title{Variational principle and almost
  quasilocality\\
  for renormalized measures.}
\author{Roberto Fern\'andez,\\
Laboratoire de Math\'ematiques Raph\"el Salem,\\
 Universit\'e de Rouen, UFR Sciences, site Colbert,\\
76821 Mont Saint-Aignan, France.\\
  E-mail: Roberto.Fernandez@univ-rouen.fr\\
  Arnaud Le Ny,\\
Laboratoire de Math\'ematiques Raph\"el Salem,\\
  Universit\'e de Rouen, UFR Sciences, site Colbert,\\
76821 Mont Saint-Aignan, France.\\
  E-mail: Arnaud.Leny@univ-rouen.fr\\
  Frank Redig,\\
  Faculteit Wiskunde en Informatica TU Eindhoven,\\
  Postbus 513, 5600 MB Eindhoven, The Netherlands.\\
  E-mail:f.h.j.redig@TUE.nl\\}

\date{July 11, 2001}

\maketitle

{\bf Keywords}: Gibbs vs non-Gibbs,  almost quasilocality.

{\bf AMS classification}:  60G60 (primary), 82B20, 82B30 (secondary).

\begin{abstract}
  We restore part of the thermodynamic formalism for some
  renormalized measures that are known to be non-Gibbsian.  We first
  point out that a recent theory due to Pfister implies that
  for block-transformed measures free energy and relative entropy
  densities exist and are conjugated convex functionals.  We then
  determine a necessary and sufficient condition for consistency with
  a specification that is quasilocal in a fixed direction. As
  corollaries we obtain consistency results for models with FKG
  monotonicity and for models with appropriate ``continuity rates''.
  For (noisy) decimations or projections of the Ising model these
  results imply almost quasilocality of the decimated ``$+$'' and
  ``$-$'' measures.
\end{abstract}

\newpage

\section{Introduction}

Non-Gibbsian measures were initially detected as ``pathologies'' of
renormalization group transformations of low-temperature Gibbs
measures \cite{GP,I}.  Initial efforts were directed towards the
construction of a sufficiently rich catalogue of examples and
mathematical mechanisms leading to non-Gibbsianness \cite{VEFS}.
the present research is focused on two practical aspects of the
phenomenon: (i) determination of weaker notions of Gibbsianness which
are preserved by the transformations of interest, and (ii) extension
of the thermodynamic formalism to these broader classes of measures.
The first issue motivated the notions of weak Gibbsianness and almost
quasilocality \cite{DOB,FP,BKL,Maes,Maes2,lef,lef2,LEN,LEN2}.
The second issue was treated in \cite{lef} and \cite{Maes2} which
discuss: (a) Existence of thermodynamic quantities, and (b) relations
between the properties ``consistency with the same specification''
(CSS) and ``zero relative entropy'' (ZRE).  Block renormalizations are
considered in \cite{lef} and projections to the line in \cite{Maes2}.
In both cases thermodynamic functionals are shown to exist, and, under
suitable hypothesis on measure supports, the implication ``CSS
$\Longrightarrow$ ZRE'' is shown in \cite{lef} and the opposed
implication in \cite{Maes2}.  Work in both references is based on the
weakly Gibbsian approach, that is on the determination of a set of
configurations of full measure on which the conditional probabilities
take a Gibbsian form, with an almost surely finite but unbounded
Hamiltonian.

Our results are complementary to those of these references on two
counts.  First, we use recent results of Pfister \cite{P} on
asymptotically decoupled measures, to restore the
specification-invariant part of the variational principle (Definition
\ref{chp:def3} below) for block-renormalized measures.  The result
(Theorem \ref{chp:thm1}) is completely independent of any
Gibbs-restoration approach.  Second, we establish conditions for the
validity of the implication ``ZRE $\Longrightarrow$ CSS'' (Theorem
\ref{r.teo1}) that, as we show in Corollary \ref{r.cor1}, can be
related to almost quasilocality.  In particular almost quasilocality
is proven for block decimations (Corollary \ref{r.cor2}),
strengthening previous weakly Gibbsian results \cite{BKL,MRSV}
obtained by more laborious means.  In Proposition \ref{prop1} we
present a further criterion that proves the implication ``ZRE
$\Longrightarrow$ CSS'' for projections to a line without further
support assumptions.

In general terms, our results illustrate the fact that important
aspects of the variational principle can be directly related to
properties of specifications, without having to rely on rather
detailed descriptions of weakly Gibbsian potentials.

\section{Basic definitions and notation}

We start by summarizing some basic notions for the sake of
completeness.  As general reference we mention \cite{HOG}.  See also
\cite{VEFS}, Section 2, for a streamlined exposition.

\subsection{Quasilocality, specifications, consistent measures}

We consider configuration spaces $\Omega=\Omega_0^{\lat}$ with
$\Omega_0$ finite and $\lat$ countable (typically $\lat=\zed$),
equipped with the product discrete topology and the product Borel
$\sigma$-algebra $\mathcal{F}$.  More generally, for (finite or
infinite) subsets $\Lambda$ of $\lat$ we consider the corresponding
measurable spaces $(\Omega_{\Lambda},\mathcal{F}_{\Lambda})$, where
$\Omega_\Lambda=\{-1,1\}^\Lambda$.  For any $\omega \in \Omega$,
$\omega_{\Lambda}$ denotes its projection on $\Omega_{\Lambda}$.  We
denote by $\mathcal{S}$ the set of the finite subsets of $\lat$.  A
function $f:\Omega\to \mathbb{R}$ is {\bf local} if there exists a
finite set $\Delta$ such that $\omega_{\Delta}=\sigma_{\Delta}$
implies $f(\omega)=f(\sigma)$.  The set of local
functions is denoted $\mathcal{F}_{\rm loc}$.  
\begin{definition}\label{def.n1}
Let $f:\Omega\to \mathbb{R}$.
\begin{itemize}
\item[(i)] $f$ is {\bf quasilocal} if it is the
uniform limit of local functions, that is, if
\begin{equation}\label{chp:eqn1}
\lim_{\Lambda \uparrow \lat} \,\sup_{\scriptstyle \sigma, \omega :\atop
\scriptstyle \sigma_{\Lambda}=\omega_{\Lambda}} \Bigl| f(\sigma)
  - f(\omega)\Bigr| \;=\;0\;.
\end{equation}
[The notation $\Lambda \uparrow \lat$ means convergence along a net
directed by inclusion.]
\item[(ii)] $f$ is {\bf quasilocal in the
  direction} $\theta\in\Omega$ if
\begin{equation}\label{chp:eqn2}
\lim_{\Lambda \uparrow \lat} \,\Bigl|
f(\omega_{\Lambda}\theta_{\Lambda^c})
- f(\omega) \Bigr| \;=\;0
\end{equation}
for each $\omega\in\Omega$.  [No uniformity required.]
\end{itemize}
\end{definition}  

Let us present some remarks.
\begin{rem}
In the present setting (product of finite single-spin spaces)
quasilocality is equivalent to continuity and uniform continuity.  This follows from
Stone-Weierstrass plus the fact that local functions are continuous
for discrete spaces.
\end{rem}

\begin{rem}
The pointwise analogues of \reff{chp:eqn1} and
\reff{chp:eqn2} are the following: (i) $f$ is \emph{quasilocal at} $\omega$ if $\lim_{\Lambda \uparrow \lat}
\,\sup_{\sigma,\eta} \Bigl| f(\omega_\Lambda\sigma_{\Lambda^c}) -
f(\omega_\Lambda\eta_{\Lambda^c}) \Bigr| =0$;\\ (ii) $f$ is
\emph{quasilocal at $\omega$ in the direction} $\theta$ if
\reff{chp:eqn2} holds only for this $\omega$.\\ We shall not resort to
these notions, but let us point out that the previous remark is no
longer valid at the pointwise level: A function can be quasilocal in
every direction at a certain $\omega$ (that is, continuous at $\omega$)
and fail to be quasilocal at $\omega$. 
\end{rem}
Here is an example illustrating the last remark.  Let $d=1$ and, for a
fixed $\omega\in\Omega$, choose a countable family $\chi^{(m)}$ of
configurations such that
$\chi^{(m)}_{[-n,n]^c}\neq \chi^{(m')}_{[-n,n]^c}$ for all
$n\in\mathbb{N}$, if $m\neq m'$, and such that $\chi^{(m)}_0\neq
\omega_0$ for all $m$.  Define
\begin{equation}
  \label{eq:nn1}
  f(\eta)\;=\; \left\{
\begin{array}{ll} 
m/(n+m) \; \hbox{if } \eta=\omega_{[-n,n]}\, \chi^{(m)}_{[-n,n]^c} 
\hbox{ for some } m,n\in\mathbb{N}\\[10pt]
0  \; \; \; \; \; \; \; \; \; \; \; \; \; \; \; \; \;\hbox{otherwise}\;.
\end{array}
\right.
\end{equation}
We see that, for all $\sigma \in \Omega$,
\begin{equation}
\lim_{\Lambda \uparrow
  \lat} f(\omega_{\Lambda}\sigma_{\Lambda^c})\to 0=f(\omega).
\end{equation}
  Hence
  $f$ is quasilocal at $\omega$ in every direction.  However,
  \begin{equation}
    \label{eq:nn3}
\sup_{\sigma,\eta} \Bigl| f(\omega_{[-n,n]}\sigma_{[-n,n]^c}) -
f(\omega_{[-n,n]}\eta_{[-n,n]^c}) \Bigr|\;=\; 
\sup_m {m\over n+m} \;=\; 1\;.
  \end{equation}
So $f$ is not quasilocal at $\omega$.
\medskip

\begin{definition}\label{rdef.sp}
A {\bf specification} on $(\Omega,\mathcal{F})$ is a family $\gamma
=\{\gamma_{\Lambda},\Lambda \in \mathcal{S}\}$ of stochastic kernels
on $(\Omega, \mathcal{F})$ that are
\begin{itemize}
\item[(I)] \label{proper} \emph{Proper}: $\forall B \in
  \mathcal{F}_{\Lambda^{c}}$,
  $\gamma_{\Lambda}(B|\omega)=\mathbf{1}_{B}(\omega)$. 
\item[(II)] \emph{Consistent}: If $\Lambda \subset \Lambda'$ are finite
  sets, then $\gamma_{\Lambda'} \gamma_{\Lambda} = \gamma_{\Lambda'}$.
\end{itemize}
\end{definition}
[We adopt the ``conditional-probability'' notation, that is,
$\gamma_{\Lambda}(A|\,\cdot\,)$ is
$\mathcal{F}_{\Lambda^{c}}$-measurable $\forall A \in \mathcal{F}$,
and $\gamma_{\Lambda}(\cdot | \omega)$ is a probability measure on
$(\Omega ,\mathcal{F})$ $\forall \omega \in \Omega$.]  The notation
$\gamma_{\Lambda'} \gamma_{\Lambda}$ refers to the natural composition
of probability kernels: $ (\gamma_{\Lambda'}
\gamma_{\Lambda})(A|\omega) =
\int_{\Omega}\gamma_{\Lambda}(A|\omega')\gamma_{\Lambda'}(d\omega' |
\omega)$.  A specification is, in fact, a strengthening of the notion
of system of proper regular conditional probabilities.  Indeed, in the
former, the consistency condition (II) is required to hold
for \emph{every} configuration $\omega\in\Omega$, and not only for almost
every $\omega\in\Omega$.  This is because the notion of specification is defined
without any reference to a particular measure.

A probability measure $\mu$ on $(\Omega,\mathcal{F})$ is said to be
\emph{consistent} with a specification $\gamma$ if the latter is a
realization of its finite-volume conditional probabilities, that is,
if $\mu[A |\mathcal{F}_{\Lambda^{c}}](\,\cdot\,)=
\gamma_{\Lambda}(A  |\,\cdot\,) \  \mu\textrm{-a.s.}$
for all $A \in \mathcal{F}$ and $\Lambda \in \mathcal{S}$.
Equivalently, $\mu$ is consistent with $\gamma$ if it satisfies the
\emph{DLR equation} (for Dobrushin, Lanford and Ruelle):
\begin{equation}
  \label{chp:eqn4}
  \mu \;=\; \mu\,\gamma_\Lambda
\end{equation}
for each $\Lambda\in \mathcal{S}$.  The right-hand side is the
composed measure: $(\mu\,\gamma_\Lambda)(f) = \int
\gamma_\Lambda(f|\omega) \,\mu(d\omega)$ for $f$ bounded measurable.
We denote $\mathcal{G}(\gamma)$ the set of measures consistent with
$\gamma$ (measures \emph{specified} by $\gamma$).  The description of
this set is, precisely, the central issue in equilibrium statistical mechanics.

A specification $\gamma$ is \emph{quasilocal} if for each $\Lambda \in
\mathcal{S}$ and each $f$ local, $\gamma_{\Lambda}f$ is a quasilocal
function.  Analogously, the specification is \emph{ quasilocal in the
  direction} $\theta$ if so are the functions $\gamma_{\Lambda}f$ for
local $f$ and finite $\Lambda$.  A probability measure $\mu$ is
\emph{quasilocal} if it is consistent with some quasilocal
specification.  Gibbsian specifications ---defined through
interactions via Boltzmann's prescription--- are the archetype of
quasilocal specifications.  Every Gibbs measure ---i.e.\ every measure
consistent with a Gibbsian specification--- is quasilocal, and the
converse requires only the additional property of non-nullness
\cite{K,Su,LEN}.  Reciprocally, a sufficient condition for
non-Gibbsianness is the existence of an \emph{essential
  non-quasilocality} (essential discontinuity), that is, of a
configuration at which \emph{every} realization of some finite-volume
conditional probability of $\mu$ is discontinuous.  These
discontinuities are related to the existence of phase transitions in
some constrained systems \cite{GP,I,VEFS}.
  
  Two categories of measures have been defined in an effort to extend
  the Gibbsian formalism to non-quasilocal measures: \emph{weakly
    Gibbsian} and \emph{almost quasilocal} (see \cite{Maes} and
  references therein).  The latter is the object of the present paper
  and we include its definition for completeness.  For a 
  specification $\gamma$ let $\Omega_{\gamma}$ be the set of
  configurations where $\gamma_{\Lambda}f$ is continuous for all
  $\Lambda \in \mathcal{S}$ and all $f$ local, and
  $\Omega_{\gamma}^\theta$ the set of configurations for which all the
  functions $\gamma_{\Lambda}f$ are quasilocal in the direction
  $\theta$.
\begin{definition}\label{chp:def4}\mbox{}
\begin{enumerate}
\item A probability measure is {\bf almost quasilocal in the
  direction $\theta$} if it is consistent with a specification $\gamma$
such that $\mu (\Omega^\theta_\gamma ) =1$.
\item A probability measure $\mu$ is {\bf almost quasilocal}
if it is consistent with a specification $\gamma$ such that
$\mu (\Omega_\gamma )=1$.
\end{enumerate}
\end{definition} 

\subsection{``Thermodynamic'' functions and the variational principle}

The variational principle links statistical mechanical and
thermodynamical quantities.  Rigorously speaking, the functions
defined below (pressure and entropy density) do not quite correspond
to standard thermodynamics.  The corresponding notions in the latter
depend only of a few parameters, while the objects below are functions
on infinite dimensional spaces.  These functions are, however, more
informative from the probabilistic point of view, because, at least in
the Gibbsian case, they are related to large-deviation principles.

Translation invariance plays an essential role in the thermodynamic
formalism.  That is, we assume that there is an action
(``translations'') $\{\tau_i:i\in\zed\}$ on $\lat$ which defines
corresponding actions on configurations ---$(\tau_i \omega)_x =
\omega_{\tau_{-i}} x$---, on functions ---$\tau_i f(\omega)=
f(\tau_{-i}\omega )$---, on measures
---$\tau_i\mu(f)=\mu(\tau_{-i}f)$--- and on specifications ---$(\tau_i
\gamma)_\Lambda(f|\omega)=\gamma_{\tau_{-i}\Lambda}
(\tau_{-i}f|\tau_{-i}\omega)$.  [To simplify the notation we will
write $\mu(f)$ instead of $E_\mu(f)$.]  Translation invariance means
invariance under all actions $\tau_i$.  We consider, in this section,
only translation-invariant probability measures on $\Omega$, whose
space we denote by $\mathcal{M}_{1,{\rm inv}}^+(\Omega)$.  We denote
$\mathcal{G}_{\rm
  inv}(\gamma)=\mathcal{G}(\gamma)\cap\mathcal{M}_{1,inv}^{+}(\Omega)$.
Furthermore, the convergence along subsets of $\lat$ is restricted to
sequences of cubes $\Lambda_n\;=\;\{\tau_{-i}(0) : i\in ([-n,n] \cap
\mathbb{Z})^d\}$.

\begin{definition}[\cite{P}]
\label{chp:def11}
A measure $\nu \in \mathcal{M}_{1,inv}^{+}(\Omega)$ is {\bf
  asymptotically decoupled} if there exist functions $g:\mathbb{N}
\longrightarrow \mathbb{N}$ and $c:\mathbb{N} \longrightarrow
[0,\infty)$ such that
\begin{equation}
\label{eq:r6}
\lim_{n \to \infty} \frac{g(n)}{n}=0 \qquad  \textrm{and} 
\qquad \lim_{n \to \infty} \frac{c(n)}{| \Lambda_n |}=0\;,
\end{equation}
such that for all $n \in \mathbb{N}$, $A \in\mathcal{F}_{\Lambda_n}$
and $B \in \mathcal{F}_{(\Lambda_{n+g(n)})^c}$,
\begin{equation}\label{chp:eqn6}
e^{-c(n)}\,\nu(A)\,\nu(B) \;\leq\; \nu(A \cap B) 
\;\leq\; e^{c(n)}\,\nu(A)\,\nu(B)\;.
\end{equation}
\end{definition}
This class of measures strictly contains the set of all Gibbs
measures.  In particular, as we observe below, it includes measures
obtained by block transformations of Gibbs measures, many of which are
known to be non-Gibbsian.

For $\mu,\nu \in (\Omega,\mathcal{F})$, the \emph{relative entropy}
at volume $\Lambda \in \mathcal{S}$ of $\mu$ relative to $\nu$ is
defined as
\begin{equation}\label{eq:rex}
H_{\Lambda}(\mu | \nu)\;=\;\left\{
\begin{array}{ll}\; \displaystyle \int_{\Omega} 
\frac{d\mu_{\Lambda}}{d\nu_{\Lambda}}
 \log\frac{d\mu_{\Lambda}}{d\nu_{\Lambda}}\, d \nu 
 \; \; \; \textrm{if} \;   \;\mu_\Lambda \ll \nu_\Lambda \\[15pt]
+ \infty \; \; \; \; \; \; \; \; \; \; \; \;\; \; \; \; \; \; \; \; \; \; \; \;\; \; \; \textrm{otherwise}.
\end{array} \right. 
\end{equation}
The notation $\mu_\Lambda$ refers to the
projection (restriction) of $\mu$ to
$(\Omega_\Lambda,\mathcal{F}_\Lambda)$.  The \emph{relative entropy
  density} of $\mu$ relative to $\nu$ is the limit
\begin{equation}\label{eq:red}
h(\mu | \nu)\;=\;\lim_{n \to \infty} \frac{H_{\Lambda_n}(\mu |
  \nu)}{| \Lambda_n |} 
\end{equation}
provided it exists.  The limit is known to exist if $\nu\in
\mathcal{M}_{1,inv}^{+}(\Omega)$ is a Gibbs measure (and $\mu\in
\mathcal{M}_{1,inv}^{+}(\Omega)$ arbitrary) and, more generally
\cite{P}, if $\nu$ is asymptotically decoupled.  In these cases
$h(\,\cdot\,|\nu)$ is an affine non-negative function on
$\mathcal{M}_{1,inv}^{+}(\Omega)$. 

For $\nu \in \mathcal{M}_{1,inv}^+(\Omega)$ and $f$ a
bounded measurable function, the \emph{pressure} (or minus free-energy
density) for $f$ relative to $\nu$ is defined as the limit
\begin{displaymath}p(f | \nu)= \lim_{n \to \infty}\,\frac{1}{|
    \Lambda_n |^{d}}\, \log\int\exp\Bigl(\sum_{x
  \in \Lambda_n}\tau_xf \Bigr)\,d\nu
\end{displaymath}
whenever it exists.  This limit exists, for every quasilocal function $f$, if
$\nu$ is Gibbsian or asymptotically decoupled \cite{P}, yielding a
convex function $p(\,\cdot\,|\nu)$.

For our purposes, it is important to separate the different
ingredients of the usual variational principle in statistical
mechanics. 
\begin{definition}[Specification-independent variational principle]
\label{chp:def3}
A measure $\nu \in \mathcal{M}_{1,inv}^+(\Omega)$ satisfies a
variational principle if the relative entropy $h(\mu | \nu)$ and the
pressure $p(f | \nu)$ exist for all $\mu \in
\mathcal{M}_{1,inv}^+(\Omega)$ and all $f \in \mathcal{F}_{{\rm loc}}$,
and they are conjugate convex functions in the sense that
\begin{equation}
  \label{eq:r2}
  p(f | \nu) \;= \; \sup_{\mu \in
  \mathcal{M}_{1,{\rm inv}}^{+}(\Omega)} \Bigl[\mu(f) - h(\mu | \nu) \Bigr]
\end{equation}
for all $f \in \mathcal{F}_{{\rm loc}}$, and
\begin{equation}
  \label{eq:r3}
  h(\mu | \nu) \;= \;\sup_{f
  \in \mathcal{F}_{{\rm loc}}}\Bigl[\mu(f) - p(f | \nu)\Bigr]
\end{equation}
for all $\mu \in \mathcal{M}_{1,inv}^{+}(\Omega)$.
\end{definition}

Gibbs measures satisfy this specification-independent principle.
Pfister \cite[Section 3.1]{P} has recently extended its validity to
asymptotically decoupled measures.  In these cases $h(\,\cdot\,|\nu)$
is the rate function for a (level 3) large-deviation principle for
$\nu$.  

\begin{definition}[Variational principle relative to a
  specification]\label{chp:defr} Let $\gamma$ be a specification and
  $\nu \in \mathcal{G}_{\rm inv}(\gamma)$.  We say that a variational
  principle occurs for $\nu$ and $\gamma$ if for all
  $\mu\in\mathcal{M}_{1,inv}^{+}(\Omega)$
\begin{equation}
\label{eq:r4}
h(\mu|\nu) =0 \ \Longleftrightarrow\ \mu\in
\mathcal{G}_{\rm inv}(\gamma)\;.
\end{equation}
\end{definition}

The equivalence \reff{eq:r4} holds for Gibbs measures $\nu$, while the
implication to the right is valid, more generally, for measures $\nu$
consistent with $\gamma$ quasilocal (see \cite{HOG}, Chapter 10).  In
\cite{lef} the implication to the left was extended to
block-transformed measures satisfying appropriate support hypothesis.
Below we extend the implication to the right to some non-Gibbsian
(non-quasilocal) measures.

\subsection{Transformations of measures}

\begin{definition}\label{r.def1}
  A {\bf renormalization transformation} $T$ from
  $(\Omega,\mathcal{F})$ to $(\Omega',\mathcal{F}')$ is a probability
  kernel $T(\,\cdot\,|\,\cdot\,)$ on $(\Omega,\mathcal{F}')$.  That
is, for each $\omega\in\Omega$ $T(\,\cdot\,|\omega)$ is a probability
measure on $(\Omega',\mathcal{F}')$ and for each $A'\in\mathcal{F}'$
$T(A'|\,\cdot\,)$ is $\mathcal{F}$-measurable.

The transformation is a {\bf block-spin transformation} if $\Omega'$
is of the form $(\Omega'_0)^{\lat'}$ and there exists $\alpha  > 0$
(compression factor) such that the following two properties hold
\begin{itemize}
\item[(i)] \emph{Strict locality}: For every $n$, $A' \in
  \mathcal{F'}_{\Lambda'_n}$ implies $T^{-1}(A') \in
  \mathcal{F}_{\Lambda'_{[\alpha n]}}$.  
\item[(ii)] \emph{Factorization}: There exists a distance ${\rm dist}$
  in $\lat'$ such that if $A' \in \mathcal{F'}_{D'}$ and $B' \in
  \mathcal{F'}_{E'}$ with ${\rm dist}(D',E') > \alpha$, then $T(A'\cap
  B'|\,\cdot\,)=T(A'|\,\cdot\,)\,T(B'|\,\cdot\,)$.
\end{itemize}
\end{definition}
A renormalization transformation is \emph{deterministic} if it is of
the form $T(\,\cdot\,|\omega)=\delta_{t(\omega)}(\,\cdot\,)$ for some
$t:\Omega\to\Omega'$.

A renormalization transformation $T$ induces a transformation $\mu
\mapsto \mu T$ on measures, with $(\mu T)(f')=\mu[T(f')]$ for each
$f'\in\mathcal{F}'$.  

In most applications, block-spin transformations have a product form:
$T(d\omega'|\omega) \;=\; \prod_{x'} T_{x'}(d\omega_{x'}|\omega)$,
where $T_x(\{\omega_{x'}\}|\,\cdot\,) \in \mathcal{F}_{B_{x'}}$, for a
family of sets $\{B_{x'}\subset\lat :x'\in\lat'\}$ ---the
\emph{blocks}--- with bounded diameter whose union covers $\lat$.
Transformations of this sort are called \emph{real-space
  renormalization transformations} in physics.  The transformations
defining cellular automata (with local rules) fit also into this
framework.  The corresponding blocks overlap and the compression
factor may be chosen arbitrarily close to one.

We briefly remind the reader of some of the transformations
considered in the sequel:
\begin{itemize}
\item \emph{Projections and decimations}: Given $D\subset\lat$, this
  is the (product) deterministic transformation defined by $t(\omega)
  = (\omega_x)_{x\in D}$.  The \emph{decimation of spacing}
  $b\in\mathbb{N}$, for which $\lat=\zed$, $D=b\zed$, is a block
  transformation, while \emph{Schonmann's example} \cite{sch89},
  corresponding to $D=$ hyperplane, is not because it fails to be
  strictly local.
\item \emph{Kadanoff}: This is a product block transformation defined
  by $T_x(d\omega_x|\omega)=\exp(p\,\omega'_x\sum_{y\in B_x}\omega_y)
  /{\rm norm}$, for a given choice of parameter $p$ and blocks $B_x$.
  In the limit $p\to\infty$ one obtains the \emph{majority
    transformation} for the given blocks.  If $B_x=bx$ this is a
  \emph{noisy decimation}, that becomes the true decimation in the
  limit $p\to\infty$.  More generally, one can define a \emph{noisy
    projection} onto $D\subset\lat$ through the transformation
  $\prod_{x\in D} \exp(p\,\omega'_x\omega_x) /{\rm norm}$.
\end{itemize}

It is well known that renormalization transformations can destroy
Gibbsianness (for reviews see \cite{VE1,VE2,fer98.1,fer98}).  Most of
the non-Gibbsian measures resulting from block transformations were
shown to be weakly Gibbsian \cite{BKL,MRSV}.  In Corollaries
\ref{r.cor1} and \ref{r.cor2} below, we show that in some instances
they are, in fact, almost quasilocal.

\subsection{Monotonicity-preserving specifications}
\label{ss.mon}

Finally we review notions related to stochastic monotonicity.  Let us
choose an appropriate (total) order ``$\le$'' for $\Omega_0$ and,
inspired by the case of the Ising model, let us call ``plus'' and
``minus'' the maximal and minimal elements.  The choice induces a
partial order on $\Omega$: $\omega\le\sigma \Longleftrightarrow
\omega_x\le\sigma_x\ \forall x\in\lat$.  Its maximal and minimal
elements are the configurations, denoted ``$+$'' and ``$-$'' in the
sequel, respectively equal to ``plus'' and to ``minus'' at each site.
For brevity, quasilocality in the ``$+$'', resp.\ ``$-$'', direction
will be called \emph{right continuity}, resp.\ \emph{left continuity}.
The partial order determines a notion of monotonicity for functions on
$\Omega$.  A specification $\pi$ is \emph{monotonicity preserving} if
for each finite $\Lambda\subset\lat$, $\pi_\Lambda f$ is increasing
whenever $f$ is.  These specifications have a number of useful
properties. In the following lemma, we summarize the properties of monotonicity preserving specifications which we need in the sequel. Proofs and more details can be found in \cite{FP}.

\begin{lemma}\label{l.mps}
Let $\gamma$ be a monotonicity-preserving specification 
\begin{itemize}
\item[(a)] The limits $\gamma^{(\pm)}_\Lambda(\,\cdot\,|\omega)
  =\lim_{S\uparrow\lat}\,\gamma_\Lambda(\,\cdot\,|\omega_S\pm_{S^c})$ exist
  and define two monoto\-nicity-preserving specifications,
  $\gamma^{(+)}$ being right continuous and $\gamma^{(-)}$ left
  continuous.  The specifications are translation-invariant if so is
  $\gamma$.  Furthermore, $\gamma^{(-)}(f)\le \gamma(f) \le
  \gamma^{(+)}(f)$ for any local increasing $f$, and the
  specifications $\gamma^{(+)}$, $\gamma^{(-)}$ and $\gamma$ are
  continuous on the set
  \begin{equation}
    \label{eq:r15}
    \Omega_\pm\;=\; \Bigl\{\omega\in\Omega:
\gamma^{(+)}(f|\omega)=\gamma^{(-)}(f|\omega)\;\forall f\in
\mathcal{F}_{{\rm loc}}, \Lambda\in\mathcal{S}\Bigr\}\;.
  \end{equation}
\item[(b)] The limits
  $\mu^\pm=\lim_{\Lambda\uparrow\lat}\,\gamma_\Lambda(\,\cdot\,|\pm)$
  exist and define two extremal measures
  $\mu^\pm\in\mathcal{G}(\gamma^{(\pm)})$ [thus $\mu^+$ is right
  continuous and $\mu^-$ left continuous] which are
  translationn-invariant if so is $\gamma$.  If $f$ is local and
  increasing, $\mu^-(f)\le\mu(f)\le\mu^+(f)$ for any
  $\mu\in\mathcal{G}(\gamma)$.
\item[(c)] For each (finite or infinite) $D\subset\lat$, the
  conditional expectations $\mu^+(f|\mathcal{F}_D)$ and
  $\mu^-(f|\mathcal{F}_D)$ can be given everywhere-defined
  monotonicity-preserving right, resp left, continuous versions.  In
  fact, these expectations come, respectively, from \emph{global
    specifications}, that is, from families of stochastic kernels
  satisfying Definition \ref{rdef.sp} also for infinite
  $\Lambda\subset\lat$.  Furthermore, $\mu^-(f|\mathcal{F}_D)\le
  \mu^+(f|\mathcal{F}_D)$ for each $f$ increasing.
\item[(d)] For each (infinite) $D\subset\lat$ there exist monotonicity
  preserving specifications $\Gamma^{(D,\pm)}$ such that the
  projections $\mu^\pm_D\in\mathcal{G}(\Gamma^{(D,\pm)})$ and
  $\Gamma^{(D,-)}_\Lambda(f)\le \Gamma^{(D,+)}_\Lambda(f)$ for each
  $f$ increasing.  [By (a) and (c) $\Gamma^{(D,+)}$ ($\Gamma^{(D,-)}$)
  can be chosen to be right (left) continuous and extended to a global
  specification on $\Omega_D$ with the same properties.]
\end{itemize}
\end{lemma}

Models satisfying the FKG property \cite{FKG} are the standard source
of monotonicity-preserving specifications.  This class of models
includes the ferromagnets with two- and one-body interactions (eg.\ 
Ising).  Item (d) of the lemma is potentially relevant for
renormalized measures because of the fact that a transformed
measure $\mu T$ can be seen as the projection on the primed variables
of the measure $\mu\times T$ on $\Omega\times \Omega'$ defined
by
\begin{equation}
(\mu\times T) (d\omega, d\omega' ) \;=\;
 T(d\omega'|\omega )\, \mu (d\omega )\;.
\end{equation}
To apply (d) of the lemma, however, one has to find a suitable
specification for this measure $\mu\times T$. If
$\mu\in\mathcal{G}(\gamma)$ and $T$ is a product transformation, a
natural candidate is the family $\gamma\otimes T$ of stochastic
kernels
\begin{equation}
  \label{eq:r21}
  (\gamma\otimes T)_{\Lambda\times\Lambda'}
  (d\omega_{\Lambda},d\omega'_{\Lambda'}|
\omega_{\Lambda^{\rm c}},\omega'_{(\Lambda')^{\rm c}}) \;=\; 
\prod_{\scriptstyle x' :\, x'\in\Lambda'\atop
\scriptstyle {\rm or}\ B_{x'}\cap\Lambda\neq\emptyset}
  T_{x'}(d\omega_{x'}|\omega_{B_{x'}}) \,
\gamma_\Lambda(d\omega_\Lambda|\omega_{\Lambda^{\rm c}})\;.
\end{equation}
\begin{definition}\label{def.r5}
A pair $(\gamma, T)$, where $\gamma$ is a
  specification and $T$ a product renormalization transformation, is a
  {\bf monotonicity-preserving pair} if the family $\gamma\otimes T$
  is a monotonicity-preserving specification.
\end{definition}
It does not seem to be so simple to construct such
monotonicity-preserving pairs.  The only examples we know of are pairs
for which $\gamma\otimes T$ is Gibbsian for a FKG interaction.  This
happens, for instance, for noisy projections (in particular noisy
decimations) of the Ising measure.

\section{Results}

The following result follows immediately from Definitions \ref{chp:def11} and
\ref{r.def1}. 

\begin{lemma}  If $\mu\in\mathcal{M}_{1}^+(\Omega)$ is asymptotically
  decoupled, then so is $\mu T$ for every block-spin transformation
  $T$. 
\end{lemma}
From the results of Pfister, we can then conclude the following:

\begin{theorem}\label{chp:thm1}
  Let $\mu\in\mathcal{M}_{1,inv}^+(\Omega)$ be asymptotically
  decoupled and $T$ be a block-spin transformation such that $\mu T$
  is translation-invariant. Then the renormalized measure $\mu T$
  satisfies the specification-independent variational principle of
  Definition \ref{chp:def3}.
\end{theorem}
In \cite[Section 3.4]{P} it is
showed that the relative entropy density $h(\,\cdot\,|\mu T)$ is the
large deviation rate function of the empirical measure $L_{\Lambda}=\sum_{x \in \Lambda} \delta_{\tau_x \sigma}$.\\

The next theorem states the criterion used in this paper to prove the
implication to the right in \reff{eq:r4} for non-quasilocal measures
$\nu$.
\begin{theorem}\label{r.teo1}
  Let $\gamma$ be a specification that is quasilocal in the direction
  $\theta\in\Omega$ and $\nu\in\mathcal{G}_{\rm inv}(\gamma)$.  For
  each $\Lambda\in\mathcal{S}$, $M\in\mathbb{N}$, $\Lambda\subset
  \Lambda_M$ and each local $f$, let $\gamma^{M,\theta}_\Lambda (f)$
  denote the function $\omega\to
  \gamma_\Lambda(f|\omega_{\Lambda_M}\,\theta_{\lat\setminus
    \Lambda_M})$.  Then, if $\mu\in\mathcal{M}_{1,{\rm
      inv}}^+(\Omega)$ is such that $h(\mu|\nu)=0$,
  \begin{equation}
    \label{eq:r.10}
    \mu\in \mathcal{G}_{\rm inv}(\gamma) \ \Longleftrightarrow\ 
\nu\Bigl[\, {d\mu_{\Lambda_M\setminus\Lambda} \over 
d\nu_{\Lambda_M\setminus\Lambda}}
\,\Bigl(\gamma^{M, \theta}_\Lambda (f) - \gamma_\Lambda (f)\Bigr) \Bigr]
\;\mathop{\longrightarrow}\limits_{M\to\infty}\;0
  \end{equation}
for all $\Lambda\in\mathcal{S}$ and $f\in\mathcal{F}_{\rm loc}$.
\end{theorem}
The right-hand-side of \reff{eq:r.10} shows that consistency requires
the concentration properties of $d\mu_{\Lambda_M\setminus\Lambda}/
d\nu_{\Lambda_M\setminus\Lambda}$ to beat asymptotic divergences due
to the lack of continuity of $\gamma_\Lambda$.  This imposes some
conditions on $\mu$ which are reminiscent of what happens for
unbounded spin-systems.  This analogy between unbounded spin systems
and non-Gibbsian measures is an early remark from Dobrushin.  Within
approaches based on potentials (weak Gibbsianness) these conditions
are defined and handled by cluster-expansion methods \cite{Maes2,lef}.
As we discuss below, in favorable cases montonicity arguments can be
used instead.  \bigskip

We present two applications of the previous theorem.  First we discuss
systems with monotonicity-preserving specifications.
\begin{corollary}\label{r.cor1}
  Consider a specification $\gamma$ that is monotonicity preserving
  and translation invariant.  Then, with the notation of Lemma
  \ref{l.mps},
\begin{itemize}
\item[(a)] $h(\mu^-|\mu^+) = 0$ implies that $\mu^- \in
  \mathcal{G}(\gamma^{(+)})$ and 
$\mu^-(\Omega_{\gamma^{(+)}})=\mu^-(\Omega_{\gamma^{(-)}})=1$
  (hence $\mu^-$ is almost quasilocal).
\item[(b)] For $\mu\in\mathcal{M}_{1,{\rm inv}}^+(\Omega)$, 
  $h(\mu|\nu^+)=0$ and $\mu (\Omega_{\pm} )=1$ implies
  $\mu\in\mathcal{G} (\gamma^{(+)} )$, and thus $\mu$ almost quasilocal.
\end{itemize}
Analogous results are valid interchanging ``$+$'' with ``$-$''.
\end{corollary}
By part (d) of Lemma \ref{l.mps}, and the comments immediately
thereafter, the preceding results apply when $\mu^\pm$ are the
projections (possibly noisy) of the ``plus'' and ``minus'' phases of
the Ising model.  More generally, they can be the renormalized
measures of the ``plus'' and ``minus'' measures of a
monotonicity-preserving specification whenever the specification and
the transformation form a monotonicity-preserving pair (Definition
\ref{def.r5}).

At low temperature, the decimations (possibly noisy) $\mu^+T$ and
$\mu^-T$ of the ``plus'' and ``minus'' phases of the Ising model are
non-Gibbsian \cite{I,VEFS}, that is, all specifications with which
these measures are consistent show essential discontinuities.  The
preceding corollary shows that, nevertheless, in these cases the
implication to the right of the variational principle \reff{eq:r4} can
be recovered up to a point.

If $\gamma$ is a quasilocal translation-invariant specification and
$T$ a block-spin transformation, then $h(\mu T|\nu T)=0$ for each
$\mu,\nu\in \mathcal{G}_{\rm inv}(\gamma)$ such that $\mu T$ and $\nu
T$ are translation invariant \cite[formula (3.28)]{VEFS}.  Hence, from
part (a) of the previous corollary we conclude the following.
\begin{corollary}\label{r.cor2}
  Let $\gamma$ be a quasilocal, monotonicity-preserving,
  translation-invariant specification, and $T$ a block-spin
  transformation that preserves translation invariance such that the
  pair $(\gamma,T)$ is monotonicity-preserving.  Let $\mu^\pm$ be the
  extremal measures for $\gamma$ [part (b) of Lemma \ref{l.mps}] and
  $\pi^{(\pm)}$ be the right(left)-continuous specifications such that
  $\mu^\pm T\in\mathcal{G}_{\rm inv}(\pi^\pm)$.
 
Then 
  \begin{equation}
    \label{eq:1}
\mu^-T \in  \mathcal{G}(\pi^+) \quad\hbox{ and }\quad
  \mu^-T(\Omega_{\pi^+})=\mu^-T(\Omega_{\pi^-})=1
  \end{equation}
  (hence $\mu^-T$ is almost quasilocal).  Analogous results are valid
  interchanging ``$+$'' with ``$-$''.
\end{corollary}

This corollary applies in particular for decimations (possibly noisy)
of the Ising model.  At low temperature, the renormalized measures
$\mu^+T$ and $\mu^-T$ are in general non-Gibbsian \cite{I,VEFS}, that
is, the specifications $\pi^+$ and $\pi^-$ show essential
discontinuities.  Nevertheless, the preceding corollary, together with
part (b) of Corollary \ref{r.cor1} shows that in these cases the
implication to the right of the variational principle \reff{eq:r4} can
be recovered, together with almost quasilocality.

Several remarks are in order.
\begin{rem}
The preceding corollary strengthens, for (noisy) decimation
transformations, the results of \cite{BKL,MRSV} where only
weak-Gibbsianness is proven.  Our argument is apparently simpler than
the renormalization and expansion-based procedures set up in these
references, but, of course, it does not give such a complete
description of the support of the decimated measures and it is only
restricted to models with monotonicity properties.
\end{rem}

\begin{rem}
For $d=2$, the corollary implies that \emph{all} the
  decimated measures of the Ising model are consistent with $\pi^+$
  and almost quasilocal.  This follows from the results of Aizenman
  \cite{aizun} and Higuchi \cite{japun} showing that $\mu^+$ and $\mu^-$ are
  the only extremal measures in $\mathcal{G}(\gamma)$.  
\end{rem}

\begin{rem}\label{r.lef}
  Lefevere proves in \cite{lef} the implication to the left in
  \reff{eq:r4}, for $\nu$ a block-transformed measure and $\mu$
  concentrated on an appropriate set $\widetilde\Omega\subset\Omega$
  of $\nu$-measure 1.
\end{rem}

Our second application of Theorem \ref{r.teo1} does not involve any
monotonicity hypothesis.  To formulate it we need some notation.  For
$\Lambda$ a fixed finite volume, and $f$ a local function, put
\begin{equation} 
\delta^\theta_{\Lambda,M}(f) \;=\;
\Bigl|\gamma^{M,\theta}_\Lambda (f)- \gamma_\la (f)\Bigr| 
\end{equation}
and introduce for $\epsilon  > 0$ the sets 
\begin{equation}
\label{rn.2}
A(\theta,\Lambda,f,\eps,M) \;=\; \{
\eta\in\Omega:\delta^\theta_{\Lambda,M}(f) > \eps \}\;.  
\end{equation}
If $\gamma$ is continuous in the direction $\theta$, then
$\delta^\theta_{\la,M}$ tends to zero as $M$ tends to infinity, and
hence for any probability measure $\mu$, $\mu
[A(\theta,\Lambda,f,\eps,M)]$ tends to zero as $M$ tends to infinity.
\begin{definition} Let $\alpha_M\uparrow\infty$ be an
increasing sequence of positive numbers.  Let $\mu $ be a probability
measure on $\Omega$.  We say that the specification $\gamma$ admits
$\alpha_M$ as a $\mu$-rate of continuity in the direction $\theta$ if
for all $\eps  > 0$, for all $f$ local, and for all $\la$:
\begin{equation}\label{conrate}
\limsup_{M\uparrow\infty} \frac{1}{\alpha_M}
\log \mu [A (\theta,\Lambda,f,\eps,M)] \; <  \;\; 0.
\end{equation}
\end{definition}
The following proposition shows that for a given rate of continuity,
the condition of Theorem \ref{r.teo1} will be satisfied if the
relative entropies tend to zero at the same rate.

\begin{proposition}\label{prop1}
  Let $\nu\in\mathcal{G}(\gamma)$ and suppose that $\alpha_M$ is a
  $\nu$-rate of $\theta$-continuity.  Suppose furthermore that
\begin{equation}\label{zerop}
\lim_{M\uparrow\infty}\frac{1}{\alpha_M} H_{\Lambda_M} (\mu|\nu ) \;=\; 0.
\end{equation}
Then $\mu\in\gee (\gamma )$.
\end{proposition}

This proposition applies, for instance, to Schonmann's example.
Indeed, if $\nu^+$ is the projection on a (one-dimensional) layer of
the low-temperature plus-phase of the two-dimensional Ising model,
then the estimates in \cite{Maes} imply that the monotone
right-continuous specification $\gamma^+$ (such that
$\nu^+\in\gee(\gamma^+ )$) admits $\alpha_M =M$ as $\nu^+$-rate of
right-continuity. Hence for this example we can conclude that for any
other measure $\mu$ on the layer, $h(\mu|\nu^+) =0$ implies
$\mu\in\gee (\gamma^+ )$.  This is a strengthening of part (b) of
Corollary \ref{r.cor1}.  We emphasize that such a $\mu$ can \emph{not}
be the projection $\nu^-$ of the minus Ising phase.  Indeed, while at
present the existence of $h(\nu^-|\nu^+)$ has not rigorously been
established, Schonmann's original argument \cite{sch89} implies that
$h(\nu^-|\nu^+) > 0$ if it exists.

\section{Proofs}

\subsection{Proof of Theorem \protect\ref{r.teo1}}

The hypothesis $h(\mu|\nu)=0$ implies, by
\reff{eq:rex}--\reff{eq:red}, that for $n$ sufficiently large the
$\mathcal{F}_{\Lambda_n}$-measurable function
$d\mu_{\Lambda_n}/d\nu_{\Lambda_n}$ exists.  Let's denote it
$g_{\Lambda_n}$.  For $f$ local and $\Lambda\in\mathcal{S}$, pick $M$
such that $\Lambda_M\supset \Lambda$ and $g_{\Lambda_M}$ exist and
write
\begin{equation}\label{eq:r22}
\mu(\gamma_\Lambda f - f) \;=\; A_M + B_M + C_M
\end{equation}
with
\begin{equation}\label{eq:r23}
  A_M\;=\; \mu\Bigl[\gamma_{\Lambda_M} (f) - \gamma^{M,
    \theta}_{\Lambda_M} (f)\Bigr]
\ ,\ B_M\;=\; \nu\Bigl[(g_{\Lambda_M}  - 
g_{\Lambda_M\setminus\Lambda})\, f\Bigr]
\end{equation}
and $C_M$ is the right-hand side in \reff{eq:r.10}.  We shall prove
that $A_M$ and $B_M$ go to zero.

Indeed, $\lim_{M\to\infty} A_M =0$ follows by dominated convergence,
because $\gamma$ is quasilocal in the direction $\theta$ and
$|\gamma^M_\Lambda (f)|\leq \| f \|_\infty$.

On the other hand, Csisz\'ar's inequality \cite{C}
%
 \begin{displaymath}
 H_{\Delta}(\mu | \nu) -  H_{\Delta'}(\mu | \nu) \;\geq\;
 \frac{1}{2}\left[\int_{\Omega}\Bigl|
 g_{\Delta}(\omega)-g_{\Delta'}(\omega) \Bigr| 
 d\nu(\omega)\right]^2. 
 \end{displaymath}
valid for $\Delta'\subset\Delta\in\mathcal{S}$, implies that
\begin{displaymath}
|B_M| \;\leq\; \sqrt{2}\, \| f \|_{\infty}\,
\Bigl[H_{\Delta}(\mu | \nu) -
H_{\Delta\setminus\Lambda}(\mu | \nu)\Bigr]
\end{displaymath}
for any $\Delta \supset \Lambda_M$.  But the hypothesis $h(\mu|\nu)=0$
implies that the difference in entropies in the right-hand side tends
to zero as $\Delta\uparrow\lat$, as shown in \cite{HOG} or \cite{P}.
Hence $B(M)\to_M 0$. 

\subsection{Proof of Corollary \protect\ref{r.cor1}}

It is enough to verify the right-hand side of \reff{eq:r.10} for
increasing local functions $f$ since linear combinations of these are
uniformly dense in the set of quasilocal functions.

\paragraph{Part (a)}
By Theorem
\ref{r.teo1} we only have to show that
\begin{equation}
  \label{eq:r24}
  C_M\;=\; \mu^+\Bigl[ g_{\Lambda_M\setminus\Lambda} 
\,\Bigl(\gamma^{M, +}_\Lambda (f) - \gamma_\Lambda (f)\Bigr) \Bigr] 
\;\mathop{\longrightarrow}\limits_{M\to\infty}0\;,
\end{equation}
where $g_D = d\mu^-_D/d\mu^+_D$ for $D\subset\lat$.
We first point out that
\begin{equation}
  \label{eq:r25}
  C_M \;\ge\; 0
\end{equation}
because $\gamma$ is monotonicity preserving, while
\begin{equation}
  \label{eq:r26}
  \mu^+\Bigl( g_{\Lambda_M\setminus\Lambda} 
\,\gamma^{M, +}_\Lambda (f)\Bigr) \;=\; 
\mu^-\Bigl(\gamma^{M, +}_\Lambda (f)\Bigr)
\end{equation}
because $\gamma^{M,+}_{\Lambda}f$ is
$\mathcal{F}_{\Lambda_M\setminus\Lambda}$-measurable.
On the other hand, 
\begin{equation}
  \label{eq:r27}
\mu^+\Bigl(g_{\Lambda_M\setminus\Lambda}\, \gamma_\Lambda (f)\Bigr) \;=\; 
\mu^+\Bigl[g_{\Lambda_M\setminus\Lambda}\, 
\mu^+ (\gamma_\Lambda (f)| \mathcal{F}_{\Lambda_M\setminus\Lambda}) \Bigr]
\end{equation}
where $\mu^+(\,\cdot\,|\mathcal{F}_D)$ are the conditional
expectations of part (c) of Lemma \ref{l.mps}.  By the last inequality
there,
\begin{equation}
  \label{eq:r28}
\mu^+\Bigl(g_{\Lambda_M\setminus\Lambda}\, \gamma_\Lambda (f)\Bigr) \;\ge\; 
\mu^+\Bigl[g_{\Lambda_M\setminus\Lambda}\, 
\mu^-(\gamma_\Lambda (f)|\mathcal{F}_{\Lambda_M\setminus\Lambda}) \Bigr]
\;=\; \mu^-\Bigl(\gamma_\Lambda (f)\Bigr)\;.
\end{equation}
because of the $\mathcal{F}_{\Lambda_M\setminus\Lambda}$-measurability
of $\mu^-(\gamma_\Lambda
(f)|\mathcal{F}_{\Lambda_M\setminus\Lambda})(\,\cdot\,)$.  From
\reff{eq:r25}, \reff{eq:r26} and \reff{eq:r28},
\begin{displaymath}
0 \;\leq\; C_M \;\leq\; \mu^-\Bigl(\bigl| \gamma^{M,+}_{\Lambda}
(f)-\gamma_{\Lambda}(f) \bigr| \Bigr)
\end{displaymath}
and \reff{eq:r24} follows from the right-continuity of $\gamma$ and
dominated convergence.
\medskip

\paragraph{Part (b)}
By monotonicity
\begin{eqnarray*}
  0 \; &\leq& \; \mu^+\Bigl[ g_{\Lambda_M\setminus\Lambda} 
\,\Bigl(\gamma^{M, +}_\Lambda (f) - \gamma_\Lambda (f)\Bigr) \Bigr] \\ 
\; &\leq& \;  \mu^+\Bigl[ g_{\Lambda_M\setminus\Lambda} \,
\Bigl(\gamma^{M, +}_\Lambda (f) - \gamma^{M,-}_\Lambda (f)\Bigr)\Bigr]
\\  
\; &=& \; \mu\Bigl(\gamma^{M, +}_\Lambda (f) - \gamma^{M,-}_\Lambda (f)\Bigr)
\end{eqnarray*}
where $g_D$ indicates the Radon-Nikodym density of $\mu_D$ with
respect to $\mu^+_D$, and the last equality follows from the
$\mathcal{F}_{\Lambda_M\setminus\Lambda}$-measurability of $\gamma^{M,
  +}_\Lambda (f) - \gamma^{M,-}_\Lambda (f)$.  The last line tends to
zero with $M$ by dominated convergence, because $\mu(\Omega_\pm)=1$.
This concludes the proof because of Theorem \ref{r.teo1}.

\subsection{Proof of Proposition \protect\ref{prop1}}

Let us fix a local function $f$, a finite set $\Lambda$ and some
$\epsilon > 0$.  We have
\begin{equation}
  \label{eq:2}
  \nu \left[
\frac{d\mu_{\la_M\setminus\la}}{d\nu_{\la_M\setminus\la}}
\left(\gamma^{M,\theta}_\la (f) -\gamma_\la (f)\right)\right]
\;\leq \;
\eps + 2 \| f \|_\infty \,\widetilde{\mu}_M (A^M_\epsilon)
\end{equation}
where $A^M_\epsilon$ denotes the set \reff{rn.2} and we abbreviated 
\begin{equation}
\widetilde{\mu}_M(A^M_\epsilon) \;=\; \nu
\left(\frac{d\mu_{\la_M\setminus\la}}{d\nu_{\la_M\setminus\la}}\,
  1_{A^M_\epsilon} \right) \;.
\end{equation}
By (\ref{conrate}) there exists $c > 0$ such that for
$M$ large enough,
\begin{equation}
\nu (A^M_\epsilon) \; \leq \; e^{-c\,\alpha_M}\;,
\end{equation}
hence, for $0 <  \;\delta  <  \;c$, and we can write the following inequalities:
\begin{eqnarray}
\widetilde{\mu}_M(A^M_\epsilon) \; &\leq&\;
 \frac{1}{\alpha_M \,\delta}\, \log \int \exp (\delta\, \alpha_M
\, 1_{A^M_\eps} )\,d\nu \;
 +\;\frac{1}{\alpha_M\,\delta} \,H(\tilde{\mu}_M |\nu)\nonumber\\
\; &\leq& \;  \frac{1}{\alpha_M \,\delta} \,\log 
\left(1+ e^{\alpha_M\,\delta}\, \nu (A^M_\eps ) \right)
\; +\; \frac{1}{\alpha_M \,\delta} \, H(\tilde{\mu}_M|\nu )\nonumber\\
\; &\leq& \; \frac{1}{\alpha_M \,\delta} \,e^{\alpha_M(\delta-c)} \;+ \;
\frac{1}{\alpha_M\, \delta}\, H(\tilde{\mu}_M|\nu )\;.
\label{eq:nn10}
\end{eqnarray}
By (\ref{zerop}), the last line tends to zero as $M\to\infty$.  By
\reff{eq:2}, and the fact that $\eps  > 0$ is
arbitrary, we conclude that condition \reff{eq:r.10} of Theorem
\ref{r.teo1} is satisfied, which implies that $\mu\in\gee(\gamma)$.

\section*{Acknowledgements}

We thank Charles-Edouard Pfister for giving us an early version of
\cite{P} and Aernout van Enter for comments and suggestions.  We are
grateful to Eurandom (RF and ALN) and The University of Groningen (RF)
for hospitality during the completion of this work.

\end{document}